\documentclass{article}

\def\C{{\mathbb{C}}}

\usepackage{amssymb}
\usepackage{epsfig}

\newtheorem{theorem}{Theorem}[section]

\newtheorem{corollary}[theorem]{Corollary}

\newcommand{\bbox}{\ \hfill\rule[-1mm]{2mm}{3.2mm}}

\title{Weierstrass integrability of differential equations
\thanks{The authors are partially supported by a MCYT/FEDER grant number
MTM2008-00694 and by a CIRIT grant number 2009SGR 381.}}

\author{{\sc Jaume Gin\'e \& Maite Grau}}

\date{Departament de Matem\`atica. Universitat de Lleida. \\
Avda. Jaume II, 69. 25001 Lleida, SPAIN. \\
{\rm E--mails:} {\tt gine@matematica.udl.cat} and {\tt
mtgrau@matematica.udl.cat}}

\begin{document}
\maketitle

\begin{abstract}
The integrability problem consists in finding the class of
functions a first integral of a given planar polynomial
differential system must belong to. We recall the characterization
of systems which admit an elementary or Liouvillian first
integral. We define {\it Weierstrass integrability} and we
determine which Weierstrass integrable systems are Liouvillian
integrable. Inside this new class of integrable systems there are
non--Liouvillian integrable systems.\\

\end{abstract}

\section{Introduction}

It is not always possible and sometimes not even advantageous to
explicitly write the solutions of a system of differential
equations in terms of elementary functions. In fact, Poincar\'e
begun the qualitative theory of differential equations to
well--understand the behavior of the solutions of a differential
system without their explicit knowledge. For Poincar\'e, it is
thus necessary to study the functions defined by the differential
equations by themselves and without bringing them back to simpler
functions. These thoughts induced Poincar\'e to tackle the study
of differential equations beyond an essentially different point of
view from his predecessors. His study provokes a conceptual change
on the understanding of differential equations. Sometimes, though,
it is possible to find elementary functions that are constants on
solution curves, that is, elementary first integrals. These first
integrals allow to occasionally deduce properties that an explicit
solution would not necessarily reveal, see for instance \cite{PS}.
This thought originated the modern integrability theory of
differential equations that tries to respond to the natural
question: When does a system of differential equations have a
first integral that can be expressed in terms of ``known
functions" an how does one find such an integral? The answer when
the ``known functions" are the elementary functions (i.e.
functions expressible in terms of exponentials, logarithms and
algebraic functions) was given in \cite{PS}, and when the ``known
functions" are the Liouvillian functions (i.e. functions that are
built up from rational functions using exponentiation,
integration, and algebraic functions) was given in \cite{S}. In
these two cases it is given the form of an integrating factor if
the system has these type of first integrals. In this paper we
extend the results presented in
\cite{PS} and \cite{S}.\\

In elementary courses on differential equations we consider
systems of the form
\begin{equation}
\dot{x}=\frac{dx}{dt}= P(x,y), \qquad  \dot{y}=\frac{dy}{dt}=
Q(x,y), \label{eq1}
\end{equation}
where $P$ and $Q$ are polynomials in $\C[x,y]$, $\C$ being the
complex numbers. Throughout this paper we will denote by $m=
\max\{\deg P, \deg Q \}$ the {\it degree} of system (\ref{eq1}).
Obviously, we can also express system (\ref{eq1}) as the
differential equation
\begin{equation}
\frac{dy}{dx}=\frac{Q(x,y)}{P(x,y)} \ . \label{eq}
\end{equation}
We learn that although we cannot always explicitly solve this
system, we are occasionally able to find first integrals, that is,
nonconstant functions $H(x,y)$, analytic on some nonempty open set
in $\C^2$, that are constant on the solution curves in this set.
To do this we consider the differential form
$Q(x,y)dx-P(x,y)dy=0$. If $\partial P / \partial x=- \partial Q /
\partial y$, then $H(x,y)=\int Q dx - P dy$ will be a first
integral. If $\partial P / \partial x \not = - \partial Q /
\partial y$, we are taught ad hoc methods to find an integrating
factor, that is a function $R(x,y)$ such that $\partial (RP) /
\partial x=- \partial (RQ) /
\partial y$. In case we can find such function $R$, $H(x,y)=\int RQ
dx - RP dy$ will be a first integral. For example, if $(\partial Q
/ \partial x+
\partial P / \partial y)/P$ is independent of $y$, then $R=\exp (
\int (\partial Q / \partial x+
\partial P / \partial y)/P dx)$ will be an integrating factor.

\section{Integrability problem}

We recall that the integrability problem consists in finding the
class of functions a first integral of a given system (\ref{eq1})
must belong to, see \cite{CGGL1}. For instance in \cite{Po},
Poincar\'e stated the problem of determining when a system
(\ref{eq1}) has a rational first integral. The works of \cite{PS}
and \cite{S} go in this direction since they give a
characterization of when a polynomial system (\ref{eq1}) has an
elementary or a Liouvillian first integral. A precise definition
of these classes of functions is given in \cite{PS,S}. An
important fact of their results is that invariant algebraic curves
and exponential factors play a distinguished role in this
characterization. Moreover, this characterization is expressed in
terms of the inverse integrating factor. Now, we state some
results related to integration of a system (\ref{eq1}) by means of
elementary and Liouvillian functions.

\begin{theorem}{\rm \cite{PS}} If system (\ref{eq1}) has an elementary
first integral, then there exists $\omega_0, \omega_1, \ldots,
\omega_n$ algebraic over the field $\C(x,y)$ and
$c_1,c_2,\ldots,c_n$ in $\C$ such that the elementary function
\begin{equation}\label{eq3}
\tilde{H}=\omega_0+ \sum_{i=1}^n c_i \ln (w_i),
\end{equation}
is a first integral of system (\ref{eq1}).
\end{theorem}

The existence of an elementary first integral is intimately
related to the existence of an algebraic inverse integrating
factor, as the following result shows.

\begin{theorem}\label{theo1}{\rm \cite{PS}} If system (\ref{eq1}) has an elementary
first integral, then there is an inverse integrating factor of the
form
\[
V= \left( \frac{A(x,y)}{B(x,y)}\right)^{1/N},
\]
where $A$, $B \in \C[x,y]$ and $N$ is a nonnegative integer
number.
\end{theorem}

In the work \cite{CGGL}, the systems (\ref{eq1}) with a
(generalized) Darboux first integral, that is, with a first
integral of the form
\begin{equation} \label{eq2}
H=f_1^{\lambda_1}f_2^{\lambda_1} \cdots f_r^{\lambda_r}
\left(\exp\left(\frac{h_1}{g_1^{n_1}}\right)\right)^{\mu_1}
\left(\exp\left(\frac{h_2}{g_2^{n_2}}\right)\right)^{\mu_2} \cdots
\left(\exp\left(\frac{h_\ell}{g_\ell^{n_\ell}}\right)\right)^{\mu_\ell},
\end{equation}
where $f_i, g_i, h_i \in \mathbb{C}[x,y]$, $\lambda_i, \mu_i \in
\mathbb{C}$ for $\forall i$ and $n_i \in \mathbb{N}$ for
$i=1,\ldots,\ell$, are studied and the following result is
accomplished.

\begin{theorem}\label{theo2}{\rm \cite{CGGL}} If system (\ref{eq1}) has a (generalized) Darboux
first integral of the form (\ref{eq2}), then there is a rational
inverse integrating factor, that is, an inverse integrating factor
of the form
\[
V= \frac{A(x,y)}{B(x,y)},
\]
where $A$, $B \in \C[x,y]$.
\end{theorem}

Unfortunately, not all the elementary functions of the form
(\ref{eq3}) are of (generalized) Darboux type. That's why we can
find systems with an elementary first integral and without a
rational inverse integrating factor. We remark that both Theorems
\ref{theo1} and \ref{theo2} give necessary conditions to have an
elementary or (generalized) Darboux, respectively, first integral.
The reciprocal to the statement of Theorem \ref{theo1} is not
necessarily true. But the reciprocal to the statement of Theorem
\ref{theo2} is true as we
will see later.\\

The following Theorem \ref{theo3} ensures that given a
(generalized) Darboux inverse integrating factor, there is a
Liouvillian first integral. The Liouvillian class of functions
contains the rational, Darboux and elementary classes of
functions. Singer gives in \cite{S} the characterization of the
existence of a Liouvillian first integral for a system (\ref{eq1})
by means of an integrating factor.

\begin{theorem} \label{theo3} {\rm \cite{S}} System (\ref{eq1})
has a Liouvillian first integral if, and only if, there is an
inverse integrating factor of the form $V=\exp \left \{
\int_{(x_0,y_0)}^{(x,y)} \eta \right \}$, where $\eta$ is a
rational $1$-form such that $d \eta \equiv 0$.
\end{theorem}

Taking into account Theorem \ref{theo3}, Christopher in \cite{Ch}
gives the following result, which makes precise the form of the
inverse integrating factor.

\begin{theorem} \label{theo4} {\rm \cite{Ch}} If the system (\ref{eq1})
has an inverse integrating factor of the form $V=\exp \left \{
\int_{(x_0,y_0)}^{(x,y)} \eta \right \}$, where $\eta$ is a
rational $1$-form such that $d \eta \equiv 0$, then there exists
an inverse integrating factor of system (\ref{eq1}) of the form
\begin{equation}
V=\exp \{D/E \} \prod C_i^{l_i}, \label{inv}
\end{equation}
where $D$, $E$, and the $C_i$ are polynomials in $x$ and $y$ and
$l_i \in \C$.
\end{theorem}

We notice that $C_i=0$ and $E=0$ are invariant algebraic curves
and $\exp \{D/E \}$ is an exponential factor for system
(\ref{eq1}), see for instance \cite{CGGL}. Theorem \ref{theo4}
states that the search of Liouvillian first integrals can be
reduced to the search of invariant algebraic curves and
exponential factors. A result to clarify the easiest functional
class of the first integral once we know the inverse integrating
factor is a straightforward consequence of Theorem \ref{theo4} and
is the reciprocal of Theorem \ref{theo2}.

\begin{corollary} \label{cor} If system (\ref{eq1}) has a
rational inverse integrating factor, then the system has a
(generalized) Darboux first integral. \end{corollary}

The proof is based in that if $V$ is a rational inverse
integrating factor of system (\ref{eq1}) then,
$\eta=(Q(x,y)dx-P(x,y)dy)/V(x,y)$ is a rational $1$-form such that
$d \eta \equiv 0$. Since $H=\exp \left \{ \int_{(x_0,y_0)}^{(x,y)}
\eta \right \}$ is a first integral of system (\ref{eq1}), by
Theorem \ref{theo4} $H$ is a (generalized) Darboux first integral.\\

In \cite{Pa}, Painlev\'e proved the following result, see also
\cite{GaGG,GG} and references therein.

\begin{theorem} \label{theo5} {\rm \cite{Pa}}
A differential system (\ref{eq1}) has a first integral of the form
\begin{equation} I(x,y) =
(y-g_1(x))^{\alpha_1}(y-g_2(x))^{\alpha_2}\ldots
(y-g_\ell(x))^{\alpha_\ell} h(x) \ , \label{trasintegral}
\end{equation}
where $g_j(x)$ are unknown particular solutions of (\ref{eq}),
$h(x)$ is an unknown function of $x$ and the $\alpha_i$ are
unknown constants such that $\prod_{i=1}^{\ell} \alpha_i \neq 0$,
if and only if it has an integrating factor of the form
\begin{equation}
M(x,y)= \frac{\alpha(x) S(x,y)}{(y-g_1(x))(y-g_2(x))\ldots
(y-g_\ell(x))} \ , \label{invers}
\end{equation}
where $S(x,y)$ is polynomial in the variable $y$ of degree
$\ell-m-1$. Moreover,
\begin{itemize}
\item[(a)] if the system has two different integrating factors
$M_1$ and $M_2$ of the form (\ref{invers}) with $M_2/M_1$
nonconstant, then there exists a change of variable that is
rational in the variable $y$ which transforms the equation
(\ref{eq}) into a Riccati equation.

\item[(b)] if the differential system has only one
integrating factor of the form (\ref{invers}), then the particular
solutions $g_i(x)$ from the ansatz (\ref{trasintegral}) are
calculated algebraically and $h(x)$ is given by a logarithmic
quadrature.
\end{itemize}
\end{theorem}

As usual we define $\mathbb{C}[[x]]$ the set of the formal power
series in the variable $x$ with coefficients in $\mathbb{C}$ and
$\mathbb{C}[y]$ the set of the polynomials in the variable $y$
with coefficients in $\mathbb{C}$. A polynomial of the form
\[
\sum_{i=0}^n a_i(x) y^i \in \mathbb{C}[[x]][y],
\]
is called a {\it formal Weierstrass polynomial} in $y$ of degree
$n$ if and only if $a_n(x)=1$ and $a_i(0)=0$ for $i<n$. A formal
Weierstrass polynomial whose coefficients are convergent is called
{\it Weierstrass polynomial}, see \cite{C}. In a natural
generalization we call {\it Weierstrass rational function} a
function which is a quotient of sums of Weierstrass polynomials.
We say that a system (\ref{eq1}) is {\it Weierstrass integrable}
if system (\ref{eq1}) admits an inverse integrating factor of the
form
\begin{equation}
V=\exp \{D/E \} \prod C_i^{l_i}, \label{inv2}
\end{equation}
where $D$, $E$, and the $C_i$ are Weierstrass polynomials and $l_i
\in \C$. In this sense, the systems which are Liouvillian
integrable are included in the set of systems which are
Weierstrass integrable because they have an inverse integrating
factor of the form (\ref{inv2}). However, there are systems which
are Weierstrass integrable which are not Liouvillian integrable,
see Example 2 in \cite{GaGG}. The systems with a first integral of
the form (\ref{trasintegral}) are Weierstrass integrable because
they have an integrating factor of the form (\ref{invers}).

The following theorem determines some Weierstrass integrable
systems which are Liouvillian integrable.

\begin{theorem} \label{theo6}
If a differential system (\ref{eq1}) has a first integral of the
form (\ref{trasintegral}), and at least one algebraic solution,
then it has a Liovillian first integral and therefore a
(generalized) Darboux inverse integrating factor of the form
(\ref{inv}).
\end{theorem}

\noindent{\it Proof}. If the system has a first integral of the
form (\ref{trasintegral}), then by Theorem \ref{theo5} either
there exists a rational Weierstrass change which transforms the
system of differential equations (\ref{eq1}) into a Riccati
equation (see statement (a)) or all the particular solutions
$g_i(x)$ from the ansatz (\ref{trasintegral}) are calculated
algebraically and $h(x)$ is given by a logarithmic quadrature (see
statement (b)). In the first case, as we have an algebraic
solution, the Riccati equation can be solved by quadratures and
the system has a Liouvillian first integral. In the second case
all the curves $y-g_i(x)=0$ are algebraic curves and $h(x)$ is
given by a logarithmic quadrature, which implies that the inverse
integrating factor (\ref{invers}) either is a rational integrating
factor if $\alpha(x)$ is a polynomial or is a (generalized)
Darboux integrating factor if $\alpha(x)$ is not a polynomial. In
both cases by Corollary \ref{cor} or Theorem
\ref{theo4}, the system has a Liouvillian first integral. \bbox \\

In the following examples we will see the existence of planar
polynomial systems which are not Weierstrass integrable.\\

In \cite{GL} it is proved that all the rational Abel differential
equation known as solvable in the literature can be reduced to a
Riccati differential equation or to a first-order linear
differential equation through a change with a rational map.
Several examples of Abel differential equations which are not
Weierstrass integrable appear in the Appendix of \cite{GL}. For
instance the Abel differential equation
\begin{equation}\label{cas2}
\frac{d y}{d x}=y^3-2xy^2,
\end{equation}
has the following first integral
\[
H(x,y)=\frac{x{\rm Ai}\left(x^2-\frac{1}{y}\right)+{\rm
Ai}\left(1,x^2-\frac{1}{y}\right)} {x{\rm
Bi}\left(x^2-\frac{1}{y}\right)+{\rm
Bi}\left(1,x^2-\frac{1}{y}\right)},
\]
where ${\rm Ai}$ and ${\rm Bi}$ are a pair of linearly independent
solutions of the Airy differential equation $\omega''=z \omega$.
However, the change $X=x^2-1/y$, $Y=x$ transforms equation
(\ref{cas2}) into the Riccati equation $dY/dX=Y^2-X$ and in these
new variables the system is Weierstrass integrable. The same
happens for all the other cases given in \cite{GL}.

In \cite{GS} it is presented a new algorithm to detect analytic
integrability or a singular expansion of the first integral around
a singular point for planar vector fields. It is also studied the
following example
\begin{equation}\label{ex} \dot{x}= -y, \quad
\dot{y}= a x + b y + y^2,
\end{equation}
and it is proved that it has a first integral of the form
\[
H(x,y)= \sum_{k=0}^{\infty} \frac{a \, e^{x}}{y^{k-1}} \, P_k(x),
\]
where $P_k$ is a polynomial of degree $\le k$. Hence, system
(\ref{ex}) is not Weierstrass integrable and for this case it is
unknown if there exists a change of variables which transforms the
system into a Weierstrass integrable one.

\end{document}